\documentclass[11pt]{article}
\usepackage{amsfonts}
\usepackage{amssymb}
\usepackage{amsmath}
\usepackage{amsthm}
\usepackage{latexsym}
\usepackage[usenames]{color}

\newtheorem{thm}{Theorem}[section]
\newtheorem{lemma}[thm]{Lemma}
\newtheorem{prop}[thm]{Proposition}
\newtheorem{cor}[thm]{Corollary}
\newtheorem{defi}[thm]{Definition}

\newtheorem*{thm*}{Theorem}
\theoremstyle{definition}
\newtheorem*{remark}{Remark}
\numberwithin{equation}{section}

\newcommand{\al}{\alpha}
\newcommand{\lp}{\left(}
\newcommand{\rp}{\right)}

\newcommand\C{{\mathbb C}}
\newcommand\N{{\mathbb N}}

\newcommand{\R}{\mathbb{R}}

\newcommand{\abs}[1]{\left\lvert #1\right\rvert}
\newcommand{\norm}[1]{\left\lVert #1\right\rVert}
\newcommand{\pd}[2]{\frac{\partial {#1}}{\partial {#2}}}

\newcommand{\Aa}{{\mathcal A}}
\newcommand {\p}   {\partial}

\def\LL{{\mathcal L}}

\usepackage{color}

\begin{document}

\title{\bf Classical Solutions for a nonlinear Fokker-Planck \\equation arising in Computational Neuroscience}

\author{
\vspace{5pt} Jos\'e A. Carrillo$^{1}$, Mar\'ia d. M. Gonz\'alez$^{2}$, Maria P. Gualdani$^{3}$ and Maria E. Schonbek$^{4}$ \\
\vspace{5pt}\small{$^{1}$ Instituci\'o Catalana de Recerca i Estudis Avan\c cats and Departament de Matem\`atiques}\\[-8pt]
\small{Universitat Aut\`onoma de Barcelona, E-08193 Bellaterra, Spain}\\
\vspace{5pt}\small{$^{2}$ ETSEIB - Departament de Matematica Aplicada I}\\[-8pt]
\small{Universitat Polit\`ecnica de Catalunya, E-08028 Barcelona, Spain}\\
\vspace{5pt}\small{$^{3}$ Department of Mathematics}\\[-8pt]
\small{UT Austin, Austin, TX 78712, USA}\\
\vspace{5pt}\small{$^{4}$ Department of Mathematics}\\[-8pt]
\small{UC Santa Cruz, Santa Cruz, CA 95064, USA}\\
}
\date{}

\maketitle

\abstract{In this paper we analyze the global existence of
classical solutions to the initial boundary-value problem for a
nonlinear parabolic equation describing the collective behavior of
an ensemble of neurons. These equations were obtained as a
diffusive approximation of the mean-field limit of a stochastic
differential equation system. The resulting Fokker-Planck equation
presents a nonlinearity in the coefficients depending on the
probability flux through the boundary. We show by an appropriate
change of variables that this parabolic equation with nonlinear
boundary conditions can be transformed into a non standard
Stefan-like free boundary problem with a source term given by a
delta function. We prove that there are global classical solutions
for inhibitory neural networks, while for excitatory networks we
give local well-posedness of classical solutions together with a
blow up criterium. Finally, we will also study the spectrum for
the linear problem corresponding to uncoupled networks and its
relation to Poincar\'e inequalities for studying their asymptotic
behavior.}


\section{Introduction}

The basic models for the collective behavior of large ensemble of
interacting neurons are based on systems of stochastic
differential equations. Each subsystem describes an individual
neuron in the network as an electric circuit model with a choice
of parameters such as the membrane potential $v$, the
conductances, the proportion of open ion channels and their type.
The individual description of each neuron includes an stochastic
current due to the spike events produced by other neurons at the
network received through the presynaptic connections. We refer to
the classical references \cite{lapicque,GK,T} and the nice brief
introduction \cite{G} for a wider overview of this area and
further references. As a result of the coupling network, the
collective behavior of the stochastic differential system can lead
to complicated dynamics: several stationary states with different
stability properties and bifurcations, synchronization, and so on,
see \cite{BrHa,NKKRC10,NKKZRC10} for instance.

To understand this behavior the evolution in time of the potential
through the cell membrane $v(t)$ has been modeled by several
authors \cite{BrHa,brunel,RBW,CBGW,sirovich,omurtag}. The neurons
relax towards their resting potential $v_L$ (leak potential) in
the absence of any interaction. All the interactions of the neuron
within the network are modeled by an incoming presynaptic current
$I(t)$ given by an stochastic process to be specified below.
Therefore, the evolution of the membrane potential is assumed to
follow the equation
\begin{equation}\label{lif}
C_m \frac{dv}{dt} = -g_L (v-v_L) + I(t) \,,
\end{equation}
where $C_m$ is the capacitance of the membrane and, $g_L$  the
leak conductance. If the voltage achieves the so-called threshold
voltage (maximum voltage), the neuron voltage is instantaneously
reset to a fixed voltage $v_R$. At each reset time the process
produces a spike, which builds up the incoming presynaptic current
$I(t)$ and is added to the mean firing rate produced by the
network $N(t)$ defined as the average number of spikes per unit
time produced in the network.

Most of the microscopic models for neuron dynamics assume that the
spike appearance times in the network follow an independent
discrete Poisson process with constant probability of emitting a
spike per unit time $\nu$. We will assume that there are two types
of neurons: inhibitory and excitatory, and each produce a spike of
strength $J_E$ and $J_I$ respectively at their spike times. The
total presynaptic current $I(t)$ in \eqref{lif}, coming from the
spikes within the network, is computed as the difference of the
total spike strengths received through the synapsis by a neuron at
the network composed by $C_E$ excitatory and $C_I$ inhibitory
neurons. This stochastic process $I(t)$ has mean given by $\mu_C=B
\nu$ with $B=C_E J_E - C_I J_I$ and, and variance $\sigma_C^2=(C_E
J_E^2+C_I J_I^2)\nu$. We will say that the network is excitatory
if $B>0$ (inhibitory respectively  if $B<0$ ). Dealing with these
discrete Poisson processes can be difficult and thus, an
approximation was proposed in the literature. This approximation
consists in substituting the stochastic process $I(t)$ by a
standard drift-diffusion process with the same mean and variance
$$
I(t)\,dt\approx \mu_C \,dt + \sigma_C\, d{\mathcal W}_t
$$
where ${\mathcal W}_t$ is the standard Brownian motion. We refer
for more details of this approximation to
\cite{BrHa,brunel,RBW,CBGW,sirovich,omurtag,mg}. The approximation
to the original Leaky Integrate\&Fire neuron model \eqref{lif} is
then given by
\begin{equation}\label{lifdef}
dv = (-v +v_L +\mu_C) \,dt + \sigma_C\, d{\mathcal W}_t
\end{equation}
where we choose the units such that $C_m=g_L=1$, for $v\leq
v_{th}$ with the jump process: $v(t_o^+)=v_R$ whenever at $t_o$
the voltage achieves the threshold value $v(t_o^-)=v_{th}$; with
$v_L<v_R<v_{th}$. The last ingredient of the model is given by the
probability of firing per unit time of the Poissonian spike train
$\nu$, i.e., the so-called total firing rate. The firing rate
depends on the activity of the network and some external stimuli,
it is  given by $\nu = \nu_{ext} + N(t)$ where $N(t)$ is the mean
firing rate produced by the network and $\nu_{ext}\geq 0$ is  the
external firing rate. The value of $N(t)$ is then computed as the
flux of neurons across the threshold or firing voltage $v_{th}$.

Studying the stochastic problem \eqref{lifdef} with the jump
process specified above can be written in terms of a partial
differential equation for the evolution of the probability density
$p(v,t)\geq 0$ of finding neurons at a voltage $v\in
(-\infty,v_{th}]$ at a time $t\geq 0$. This PDE has the structure
of a backward Kolmogorov or Fokker-Planck equation with sources
and given by
\begin{equation}
\label{eq:nif1} \frac{\partial p}{\partial t} (v,t) =
\frac{\partial}{\partial v} \left[\big(v -v_L - \mu_C\big)
p(v,t)\right] + \frac{\sigma^2_C}{2} \frac{\p^2 p}{\p  v^2} (v,t)
+ N(t) \, \delta_{v=v_R}, \qquad v\leq v_{th} \, .
\end{equation}
A Delta Dirac source term in the right-hand side appears due to
the  firing  at time $t\geq 0$  for neurons whose voltage is
immediately reset to $v_R$. Imposing the condition that no neuron
should have the firing voltage due to their instantaneous
discharge, we complement \eqref{eq:nif1} with Dirichlet and
initial boundary conditions
\begin{equation} \label{eq:nif2}
p(v_{th},t)=0, \qquad p (-\infty,t)=0, \qquad p(v,0)=p_I(v)\geq
0\, .
\end{equation}
The mean firing rate $N(t)$ is implicitly given by
\begin{equation} \label{eq:nif3}
N(t):= -\frac{\sigma^2_C}{2} \frac{\p p}{\p v} (v_{th},t) \geq 0
\,,
\end{equation}
that is the flux of probability  of neuron's voltage,  is the least at $v_{th}$.  It is easy to check that this
definition implies, at least formally, that the evolution of
\eqref{eq:nif1} is a probability density for all times, that is
$$
\int_{-\infty}^{v_{th}} p(v,t)\,dv = \int_{-\infty}^{v_{th}}
p_I(v)\,dv =1
$$
for all $t\geq 0$. Let us note that in most of the computational
neuroscience literature \cite{BrHa,mg}, equation \eqref{eq:nif1}
is specified on the intervals $(-\infty,v_R)$ or $(v_R,v_{th})$ with
no source term but rather a boundary condition relating the values
of the fluxes from the right and the left at $v=v_R$. The
formulation presented here is equivalent and more suitable for
mathematical treatment. Other more complicated microscopic models
including the conductance and leading to kinetic-like
Fokker-Planck equations have  been studied recently, see
\cite{CCTa} and the references therein.

Finally, the nonlinear Fokker-Planck equation can be rewritten as
\begin{equation*}
\frac{\partial p}{\partial t} = \frac{\sigma^2}{2} \frac{\partial
^2p}{\partial v^2}+\frac{\partial }{\partial v} [(v-\bar\mu)P]
+N(t) \,\delta_{v=v_R},\;  v\leq v_{th}
\end{equation*}
where $\sigma^2 = 2a_0^2 +a_1N(t) $, with $a_0>0$, $a_1\geq 0$ and
$\bar\mu = B \nu_{ext}+BN(t)$. We will focus only on the simplest
case in which  the nonlinearity in the diffusion coefficient is
neglected by assuming $a_1=0$. Without loss of generality, we can
choose a new voltage variable $\tilde v\leq 0$ and an scaled
density $\tilde p$ defined by
$$
\tilde p (t, \tilde v)=\beta p (t, \beta\tilde v + v_{th})
$$
where $\beta=a_0$. Then our  main equation, after dropping the
tildes, reads
\begin{equation} \label{eq:C}
\pd{p}{t}=\frac{\partial^2 p}{\partial
v^2}+\frac{\partial}{\partial
v}\left[(v-\mu)p\right]+N(t)\,\delta_{v= v_R}, \quad \;v\leq 0\,,
\end{equation}
where the drift term, source of the nonlinearity, is given by
\begin{equation}\label{drift}
\mu = b_0+bN(t)\,\qquad \mbox{with } N(t) =-\frac{\p p}{\p v}
(0,t) \geq 0
\end{equation}
with $b_0=(B\nu_{ext}-v_{th})/a_0$ and $b=B/a_0^3$. Let us remark
that the sign of $b_0$ determines if the neurons due only to
external stimuli may produce a spike or not, therefore it controls
the strength of the external stimuli.

In a recent work \cite{Caceres-Carrillo-Perthame}, it was shown
that the problem \eqref{eq:C}-\eqref{eq:nif2}-\eqref{drift} can
lead to finite-time blow up of solutions for excitatory networks
$b>0$ and for initial data concentrated close enough to the
threshold voltage. Here, we give a characterization of the maximal
time of existence of the classical solution, if finite, and thus,
of the blow up time. We  show that if the maxima existence time is
finite, it coincides with the time in which the firing rate $N(t)$
diverges. This divergence in finite time of the firing rate has no
clear biological significance. It could mean that some sort of
synchronization of the whole network happens, see
\cite{Caceres-Carrillo-Perthame} for a deeper discussion. This is
an scenario that does not show up in the  typical reported
applications \cite{BrHa,brunel}. In the rest of this work, we
concentrate in studying the existence of classical solutions to
the initial boundary value problem
\eqref{eq:C}-\eqref{eq:nif2}-\eqref{drift}. We  show that
solutions exists globally in time for inhibitory networks $b<0$
and, we  give a characterization of the blow up time for the case
when $b>0$. Although the precise notion of classical solution will
be discussed in the next section, the main theorem of this work
can be summarized as follows.

\begin{thm}\label{main}
Let  $p_I(x) $ be a non-negative $\mathcal
C^1((-\infty,v_{th}])$ function such that $p_I(v_{th})=0$. Suppose  that
$p_I, (p_I)_x$ decay at $-\infty$, then there exists a unique
classical solution to the problem
\eqref{eq:C}-\eqref{eq:nif2}-\eqref{drift} on the time interval
$[0,T^*)$ with $T^*=\infty$ for $b\leq 0$ and, for $b>0$,  $T^*>0$
can be characterized by
$$
T^*=\sup\{t>0\,:\, N(t)<\infty\}\,.
$$
Furthermore, for $b>0$ there exist classical solutions blowing up
in finite time, and thus, with diverging mean firing rate in
finite time.
\end{thm}

Let us remark that the last statement is merely obtained by
combining the result in \cite{Caceres-Carrillo-Perthame} with our
classical solutions existence result and the characterization of
the maximal time of existence.

The main strategy of the proof as is shown in section 2, is given
by an equivalence. This equivalence, through an explicit
time-space change of variables, transforms our problem  into a
Stefan-like free boundary problem with Delta Dirac source terms,
resembling price-formation models studied in
\cite{price-formation}. In section 3, we will use ideas and
arguments in Stefan problems
\cite{Friedman:book,Friedman:free-boundary-I} to show local
existence of a solution. Next, in section 4 we will prove global
existence of classical solutions for inhibitory networks $(b<0)$
and give a characterization of the blow up time for excitatory
networks $(b>0)$. The difference between the cases $b<0$ and $b>0$
corresponds to the well studied Stefan-problem in the normal and
in the undercooled cases, see \cite{Meirmanov} for classical
references in the Stefan problem. The final section is devoted to
study the spectrum of the linear version of (\ref{eq:C}) ($b=0$)
that has some interesting features and properties connected to
classical Fokker-Planck equations.


\section{Relation to the Stefan problem}

The main aim of this section  is to rewrite equation \eqref{eq:C}
as a free boundary Stefan problem with a nonstandard right hand
side.  For this we recall a well known  change of variables,
\cite{CarrilloToscani00}, that transforms Fokker-Planck type
equations into a non-homogeneous heat equation. This change of
variables is given by
$$
y=e^tv,\;\;\tau =\frac{1}{2}(e^{2t}-1),
$$
that yields
\[
p(v,t)= e^tw\lp e^t v,\frac{1}{2}(e^{2t}-1)\rp,
\]
or equivalently
\[
w(y,\tau) =(2\tau+1)^{-1/2}p\left(\frac{y}{\sqrt{2\tau+1}}, \frac{1}{2} \log(2\tau +1)\right).
\]
 In the sequel, to simplify the notation, we use $\al(\tau) =
(2\tau +1)^{-1/2}= e^{-t}$. A straightforward computation gives
for $w$:
\begin{equation} \label{eq:w}
w_{\tau} = w_{yy} -\mu(\tau) \al(\tau) w_y +  M(\tau)\delta_{y=
\frac{v_R}{\al{(\tau)}}}
\end{equation}
where $ M(\tau) = \al^2(\tau) N(t) = -\left.\frac{\partial
w}{\partial y} \right|_{y=0}$. The additional change
of variables:
\[
u(x,\tau) =w(y,\tau)\;\; \mbox{where}\, x= y-\int_0^{\tau} \mu(s)
\al(s) \,ds = y-b_0 \left(\sqrt{1+2t}-1\right)-b\int _0^{\tau}M(s)
\alpha^{-1}(s)\,ds,
\]
removes the term with $w_y$ in \eqref{eq:w}. For boundary conditions at initial time, denote $s_I=v_{th}(=0)$. We have the following equivalent equation

\begin{lemma}
System \eqref{eq:C}-\eqref{eq:nif2}-\eqref{eq:nif3} is equivalent
to the following problem
\begin{equation}\left\{\begin{aligned} \label{stefan}
u_t &=u_{xx} + M(t)\delta_{x=s_1(t)},&  x<s(t), t>0,\\
s(t) &=s_I-b_0 \left(\sqrt{1+2t}-1\right)-b\int_0^{t} M(s) \al^{-1}(s)\,ds,& t>0,\\
s_1(t) &= s(t) +  \frac{v_R}{\al{(t)}},& t>0,\\
M(t) &=- \left.\frac{\partial u}{\partial{x}}\right|_{x=s(t)},& t>0,\\
u(-\infty,t) &=0,\;\;u(s(t),t) =0,& t>0, \\
u(x,0)&=u_{I}(x),& x<s_I.
\end{aligned}\right.\end{equation}
\end{lemma}
\begin{proof}
The proof is straightforward by the changes of variables
specified above and, as such is omitted.
\end{proof}

We now give a definition the concept of classical solution. In what follows  we work with
the Stefan-like free boundary problem \eqref{stefan}. It is
immediate to translate this to a concept of classical solution to
the original problem \eqref{eq:C}-\eqref{eq:nif2}-\eqref{eq:nif3}
by substituting $u$ by $p$, $x$ by $v$, $M(t)$ by $N(t)$, $s_1(t)$
by $v_R$, and $s(t)$ by $v_{th}$.

\begin{defi}\label{defi-solution}
Let  $u_I(x) $ be a non-negative $\mathcal C^1((-\infty,s_I])$
function such that $u_I(s_I)=0$. Suppose that $u_I, (u_I)_x$ decay at
$-\infty$. We say that $(u(x,t),s(t))$ is a solution of
\eqref{stefan} with initial data $u_I(x)$ on the time interval
$J=[0,T)$ or $J=[0,T]$, for a given $0<T\leq\infty$, if:
  \begin{enumerate}
  \item $M(t)$ is a continuous function for all $t\in J$,
  \item $u$ is continuous in the region $\{(x,t): -\infty<x\leq s(t), t\in J\}$,
   \item $u_{xx}$ and $u_{t}$ are continuous in the region $\{(x,t): -\infty<x<s_1(t), t\in J\backslash\{0\}\}\cup\{(x,t):s_1(t)<x<s(t), t\in J\backslash\{0\}\}$,
  \item $u_x(s_1(t)^-,t)$, $u_x(s_1(t)^+,t)$, $u_x(s(t)^-,t)$ are well defined,
  \item $u_x$ decays at $-\infty$,
  \item Equations \eqref{stefan} are satisfied.
 \end{enumerate}
\end{defi}
\medskip

The next lemma presents some of the a priori properties of the
solution to \eqref{stefan}.

\begin{lemma} \label{lemma-properties}
Let $u(x,t) $ be a solution to \eqref{stefan} in the sense of
Definition {\rm\ref{defi-solution}}.Then
\begin{enumerate}
\item[i)] The mass is conserved,
$$
\int_{-\infty}^{s(t)} u(x,t)\,dx= \int_{-\infty}^{s_I} u_I(x)dx,
$$
for all $t>0$.

\item[ii)] The flux across the free boundary $s_1$ is exactly the
strength of the source term:
\[
M(t):=-u_x(s(t),t)= u_x(s_1(t)^{-},t) - u_x(s_1(t)^{+},t).
\]

\item[iii)] If $b_0<0$ and $b<0$ (resp. $b_0>0$ and $b>0$), the
free boundary $s(t)$ is a monotone increasing (resp. decreasing)
function of time.
\end{enumerate}
\end{lemma}

\begin{proof} {\it i)} Mass conservation,  follows by  integration of the
equation and straightforward integration by parts. \\
{\it ii)} To establish the jump across the free boundary, i.e.
part \emph{ii)}, integrate the first equation in (\ref{stefan})
over the interval  $(-\infty, s_1(t))$ , yielding
\[
\int_{-\infty}^{s_1(t)} u_t dx-  \int_{-\infty}^{s_1(t)} u_{xx} dx =0.
\]
Hence,
\begin{equation} \label{cosa1}
\frac{\partial}{\partial t} \int_{-\infty}^{s_1(t)} u(x,t) dx =
 u_x(s_1(t)^{-},t)+\dot{s}_1(t)u(s_1(t),t).
\end{equation}
Similarly, an integration of the first equation in  \eqref{stefan}
in the interval $(s_1(t),s(t))$ gives
\begin{equation*}
\frac{\partial}{\partial t} \int_{s_1(t)}^{s(t)} u(x,t) dx +
\dot{s}_1(t) u(s_1(t),t) - \dot{s}(t) u(s(t),t)= u_x(s(t),t)
-u_x(s_1(t)^{+},t).
\end{equation*}
If we substitute $u(s(t),t) =0$ in the previous line it follows
\begin{equation} \label{cosa2}
\frac{\partial}{\partial t} \int_{s_1(t)}^{s(t)} u dx +
\dot{s}_1(t) u(s_1(t),t)= u_x(s(t),t) -u_x(s_1(t)^+,t) .
\end{equation}
Adding \eqref{cosa1} to \eqref{cosa2}  and recalling that the mass
is  preserved we get
\[
0 = \frac{\partial}{\partial t} \int_{-\infty}^{s(t)} u(x,t) dx =  u_x(s_1(t)^{-},t) + u_x(s(t),t) - u_x(s_1(t)^{+},t).
\]
It follows  that
\[
u_x(s(t),t)=   u_x(s_1(t)^{+},t) - u_x(s_1(t)^{-},t),
\]
as desired.

Let us prove the last part \emph{iii)}. The free boundary is a
monotone increasing function of time since $b_0<0$, $b<0$, $\alpha
>0$, and
\[
s(t) =s_I-b_0 \left(\sqrt{1+2t}-1\right)-b\int_0^{t} M(s)
\al^{-1}(s)\,ds, \quad t>0,
\]
while the fact that $M(t)$ is strictly positive follows by the
classical Hopf's lemma.
\end{proof}


\section{Local existence and uniqueness}

In this section we prove local existence of solution. Our method
is inspired by the theory developed by Friedman in
\cite{Friedman:book,Friedman:free-boundary-I} for the Stefan
problem. We first derive an integral formulation for the problem.
A derivative with respect to $x$  yields an integral equation for
the flux $M$, where a  fixed point argument can be used to obtain short
time  existence.  Once $M$ is known the equation for $u$ decouples
and it is solved as a linear equation.

\begin{thm} \label{short-existence} Let $u_I(x) $ be a
non-negative $\mathcal C^1((-\infty,s_I])$ function such that
$u_I(s_I)=0$, and $u_I, (u_I)_x$ decay at $-\infty$. Then there
exists a time  $T>0$, and a unique solution $(u(x,t),s(t))$ of the
equation \eqref{stefan} in the sense of Definition
{\rm\ref{defi-solution}} for $t \in [0,T]$ with initial data
$u_{I}$. Moreover, the existence time $T$ is an inversely
proportional function of
$$
\sup_{-\infty<x\leq s_I}|u_{I}'(x)| \,.
$$
\end{thm}

The proof of Theorem \ref{short-existence} will be divided in
several steps. The first step deals  with an integral formulation
of the solution, which is used to show the existence of $M$.

\subsection{The integral formulation}
Let $G$ be the the Green's function for the heat equation on the
real line:
\begin{equation*}
G(x,t,\xi,\tau)= \frac{1}{[4\pi(t-\tau)]^{1/2}}
\exp\left\{-\frac{|x-\xi|^2}{4(t-\tau)}\right\}.
\end{equation*}
To obtain an integral formulation of the solution $u$ of
\eqref{stefan}, recall the following Green's identity
\begin{equation}\label{Green-identity}
\frac{\partial}{\partial \xi} \left( G \frac{\partial u}{\partial
\xi} - u \frac{\partial G}{\partial \xi} \right)
-\frac{\partial}{\partial \tau} (Gu) =0.
\end{equation}
To recover $u$ we first integrate the identity
\eqref{Green-identity} in  the two regions
$$
-\infty  <\xi<s_1(\tau),\;0<\tau<t,\quad\mbox{and} \quad s_1(\tau) <\xi<s(\tau),\;0<\tau<t,
$$
and then add up the results from the integration. We split the
resulting expression into the following four terms; the
only problematic one is the one containing $u_{\xi\xi}$:
\begin{align*}
I &=\int_0^t \int_{-\infty}^{s_1(\tau)} \frac{\partial}{\partial
\xi}\lp G \frac{\partial u}{\partial \xi}\rp d\xi d\tau,\qquad
II =\int_0^t \int_{s_1(\tau)}^{s(\tau)} \frac{\partial}{\partial \xi}\lp G \frac{\partial u}{\partial \xi}\rp d\xi d\tau,\nonumber\\
III &= \int_0^t \int_{-\infty}^{s(\tau)} \frac{\partial}{\partial
\xi}\lp u \frac{\partial G}{\partial \xi}\rp d\xi d\tau,\qquad
IV=\int_0^t \int_{-\infty}^{s(\tau)} \frac{\partial}{\partial
\tau}(G u) \, d\xi d\tau.
\end{align*}
Each term will be analyzed separately. Note that $u$ and $G$ have
enough decay as $ |\xi | \to  \infty$ to justify the following computations
due to Definition \ref{defi-solution}. Since
$G(x,t,-\infty,\tau)=0$ it holds
\begin{equation} \label{I}\begin{aligned}
I &= \int_0^t \left. G \frac{\partial u}{\partial \xi}\right|_{\xi=-\infty}^{\xi = s_1(\tau)} d\tau
=\int_0^t G(x,t,s_1(\tau),\tau) \left. \frac{\partial u}{\partial \xi}\right|_{s_1(\tau)^{-}} d\tau.
\end{aligned}\end{equation}
Next, we obtain
\begin{equation}\label{II}
II = \int_0^t \!\! \left\{ \left.G\frac{\partial u}{\partial
\xi}\right|_{\xi=s(\tau)} \!\!\!\!\!\!\!\!- G\left.\frac{\partial
u}{\partial \xi}\right|_{\xi=s_1(\tau)^+}\right\}d\tau = -
\int_0^t \!\! \left\{ \left.G\right|_{\xi=s(\tau)} M(\tau) +
G\left.\frac{\partial u}{\partial
\xi}\right|_{\xi=s_1(\tau)^{+}}\right\}d\tau.
\end{equation}
Here we have used that $\left.\frac{\partial u}{\partial
\xi}\right|_{\xi=s(\tau)}= -M(\tau)$. For the third integral we
have
\begin{equation}\label{III}\begin{aligned}
III &= -\int_0^t \left\{ \left.\lp u \frac{\partial G}{\partial \xi}\rp\right|_{\xi =s(\tau)}-\left.\lp u \frac{\partial G}{\partial \xi}\rp\right|_{\xi=-\infty}\right\} d\tau\\
&= -\int_0^t \left\{ \left.(u(s(\tau),\tau)  \frac{\partial G}{\partial \xi}\right|_{\xi =s(\tau)} -\left.u(-\infty,\tau) \frac{\partial G}{\partial \xi}\right|_{\xi=-\infty}\right\} d\tau=0,
\end{aligned}\end{equation}
taking into account that $u(s(\tau),\tau) = u(-\infty,\tau)=0$.
Finally, using that $u(s(\tau),\tau) =0$, we have
$$
IV= \int_0^{t} \frac{\partial}{\partial \tau} \int_{-\infty}^{s(\tau)} Gu d\xi d\tau= \int_{-\infty}^{s(t)} \left.Gu\right|_{\tau=t} d\xi -  \int_{-\infty}^{s(0)} \left.Gu\right|_{\tau=0} d\xi.
$$
Recall  that $G(x,t,\xi,t) = \delta_{x=\xi}$, thus  the last
identity yields
\begin{equation} \label{IV}
IV = \int_{-\infty}^{s(t)}\delta_{\xi=x} u(\xi,t)d\xi -
\int_{-\infty}^{s(0)} G(x,t,\xi,0) u_I(\xi)  d\xi.
\end{equation}
Combining (\ref{I}), (\ref{II}), (\ref{III}), \eqref{IV}, and  part \emph{ii)} of Lemma \ref{lemma-properties}, yields that the
solution $u$ reads as
\begin{align}\label{Duhamel}
u(x,t) =& \int_{-\infty}^{s(0)} G(x,t,\xi,0) u_I(\xi)  d\xi + \int_0^t G(x,t,s_1(\tau) )\left.\frac{\partial u}{\partial \xi}\right|_{\xi=s_1(\tau)^{-}}\nonumber\\
&- \int_0^t M(\tau) G(x,t,s(\tau),\tau) d\tau - \int_0^t \left.G(x,t,s_1(\tau) )\frac{\partial u}{\partial \xi}\right|_{\xi=s_1(\tau)^{+}}\nonumber\\
=&\int_{-\infty}^{s(0)} G(x,t,\xi,0) u_I(\xi)  d\xi -\int_0^t M(\tau) G(x,t,s(\tau),\tau) d\tau + \int_0^t M(\tau)G(x,t,s_1(\tau),\tau)d\tau\nonumber\\
 =&:\;  I_1-I_2+I_3.
\end{align}
 The term $I_1$
represents the solution to the heat equation with data $u_I$.
Indeed,
\[\int_{-\infty}^{s(0)} G(x,t,\xi,0) u_{I}(\xi) d\xi \] is the solution to the homogeneous heat equation with initial data
\[u_0(\xi) =
\left\{
\begin{array}{lr}
u_I(\xi)&\xi \le s(0),\\
0&\xi > s(0).
\end{array}
\right.
\]

All the calculations up to here are formal assuming that $u$ is a
solution of the equation \eqref{stefan} as in Definition
\ref{defi-solution}. We now derive an equation for $M$ which will
be solved for short time using a fixed point argument. The first
step is to obtain the space derivatives of the terms
$I_i,\;i=1,2,3$ and evaluate them at $x=s(t)^{-}$:
\begin{align*}
\left.\frac{\partial I_1}{\partial x} \right|_{x=s(t)^{-}} = \int_{-\infty}^{s(0)} G_x(x,t,\xi,0) u_I(\xi)  d\xi=- \int_{-\infty}^{s(0)}G(x,t,\xi,0) u'_I(\xi)  d\xi .\notag
\end{align*}
To get the derivative of  $I_2$, we use \cite[Lemma 1, pag
217]{Friedman:book}: this lemma states  that for any continuous
function $\rho$,
\begin{equation}\label{little-lemma}
\lim_{x\to s(t)^-} \frac{\partial}{\partial x} \int_0^t \rho(\tau)
G(x,t,s(\tau),\tau)\,d\tau=\frac{1}{2} \rho(t)+\int_0^t \rho(\tau)
\frac{\partial G}{\partial x} (s(t),t,s(\tau),\tau)\, d\tau.
\end{equation}
As a consequence,
\[
\left.\frac{\partial I_2}{\partial x} \right|_{x=s(t)^{-}} =  \frac{1}{2} M(t) +\int_{0}^t M(\tau) G_x(s(t),t,s(\tau),\tau) d\tau.
\]
For  the derivative of $I_3$  note that problems can only occur if
$t=\tau$ and $s(t) =s_1(\tau)$, but this is not possible  by the definition of $s_1$
. Thus,
\[
\left.\frac{\partial I_3}{\partial x} \right|_{x=s(t)^{-}} = \int_0^t G_x(s(t),t;s_1(\tau),\tau) M(\tau) d\tau.
\]
Substituting the  estimates on $I_1$, $I_2$ and $I_3$
into (\ref{Duhamel}) we get
\begin{align*}
-M(t) =& \int_{-\infty}^{s(0)}G(s(t),t,\xi,0) u_I'(\xi)  d\xi -\frac{1}{2} M(t)\\
&-\int_{0}^t M(\tau)G_x(s(t),t,s(\tau),\tau) d\tau +\int_0^t M(\tau) G_x(s(t),t;s_1(\tau),\tau)  d\tau.
\end{align*}
After reordering yields,
\begin{equation}\label{M}\begin{aligned}
 M(t) =& -2\int_{-\infty}^{s(0)}G(s(t),t,\xi,0) u'_I(\xi) \, d\xi \\
 &+2\int_{0}^t M(\tau)G_x(s(t),t,s(\tau),\tau)\, d\tau - 2\int_0^t M(\tau) G_x(s(t),t,s_1(\tau),\tau) \, d\tau.
\end{aligned}\end{equation}


\subsection{Local existence  and uniqueness for $M$}

\begin{thm} \label{thm:local-M}
Let $u_I(x) $ be a non-negative $\mathcal C^1((-\infty,s_I])$
function such that $u_I(s_I)=0$. Suppose  $u_I,(u_I)_x$ decay to zero
as $ x\to -\infty$. Then there exists a time $T>0$ such that
$M(t)$ defined by the integral formulation \eqref{M} exists for
$t\in [0,T]$ and is unique in $\mathcal C([0,T])$. The
existence time $T$ satisfies
$$
T \leq \left(\sup_{-\infty<x\leq s_I}|u_{I}'(x)| \right)^{-1}\,.
$$
\end{thm}

\begin{proof} The local in time existence of $M(t)$ is obtained via a
fixed point argument.  For this, we  modify the classical
argument for the Stefan problem to account for the additional
source term given by $M(t)\,\delta_{x=s_1(t)}$. For $\sigma,m>0$,
consider the norm
$$
\norm{M}:= \sup_{0\leq t\leq\sigma} |M(t)|
$$
in the space
\[
C_{\sigma,m} :=\{M\in \mathcal C([0,\sigma]) : \|M\| \leq m\}.
\]
Set
\begin{equation}\label{T}\begin{aligned}
T(M)(t) :=&  -2\int_{-\infty}^{s(0)}G(s(t),t,\xi,0) u'_I(\xi)  d\xi \\
& +2\int_{0}^t M(\tau)G_x(s(t),t,s(\tau),\tau) d\tau - 2\int_0^t M(\tau) G_x(s(t),t,s_1(\tau),\tau)  d\tau\\
:=& \; J_1 +J_2+J_3.
\end{aligned}\end{equation}
In order to apply fixed point arguments, it is necessary to show
that for  sufficiently small $\sigma$ we have: $T:C_{\sigma,m}\to
C_{\sigma,m}$ and that $T$ is a contraction. Define
\begin{equation}\label{choice-m}
m:= 1 +2 \sup_{-\infty<x\leq s(0)}|u_{I}'(x)|.
\end{equation}

{\it Step 1.-} We show that for $\sigma$ sufficiently small  $T:C_{\sigma,m}\to C_{\sigma,m}$. For simplicity, we focus on  the proof in
the case $b<0$. At the end we make the necessary corrections
for $b>0$. Choose $\sigma$ sufficiently small so that
\begin{itemize}
\item[\emph{i.}] $\al^{-1}(t) \leq 2,\;\forall \; t\leq \sigma$ ,

\item[\emph{ii.}]$ \frac{m(|b_0|+2m|b|)}{\sqrt{\pi} }\sigma^{1/2}
\leq 1/2$,

\item[\emph{iii.}]$|v_R|-|b_0|\sigma >0$,

\item[\emph{iv.}] $\displaystyle
\frac{2m}{\sqrt{\pi}}\int_{\frac{\abs{v_R}-|b_0|\sigma}{\sqrt{8\sigma}}}^{\infty}
z^{-1}\exp\{-z^2\} dz \leq 1/2.$
\end{itemize}
We obtain first an auxiliary estimate. Since $\sigma$ has been chosen so small that
condition \emph{i.} holds and $\alpha^{-1}\sqrt{1+2t}$ is a 1-Lipschitz
function for $t\geq 0$, if $M\in C_{\sigma,m}$ then
\begin{equation}\label{s-Lipschitz}
|s(t)-s(\tau)| \leq |b_0||t-\tau|+ |b|\int_{\tau}^t M(s)
\al^{-1}(s)\,ds \leq \left(|b_0|+2|b|m\right) |t-\tau|,
\end{equation}
i.e., $s(t)$ is a Lipschitz continuous function of time.

To estimate the image of the operator $T(M)$ as defined in
\eqref{T} for $M\in  C_{\sigma,m}$ we find separately a bound for
each term $J_1, J_2, J_3$. First, for $J_1$, note that
$$
\int_{-\infty}^{s(0)} G(s(t),t,\xi,0) \,d\xi \leq 1.
$$
Then, it is straightforward to check
\[
\abs{J_1} \leq 2\left\{\sup_{-\infty<x\leq s(0)}|u_{I}'(x)|\right\}\int_{-\infty}^{s(0)} G(x,t,\xi,0) d\xi \leq 2\sup_{-\infty<x\leq s(0)}|u_{I}(x)|.
\]
We bound $J_2$ as
$$
\abs{J_2} \leq 2m\int_0^t |G_x(s(t),t,s(\tau),\tau)| d\tau.
$$
Substituting
$$
G_x (x,t,\xi,\tau)= -\frac{1}{2\sqrt{4\pi}}
\frac{(x-\xi)}{(t-\tau)^{3/2}}\exp\left\{-\frac{|x-\xi|^2}{4(t-\tau)}\right\}
$$
in \eqref{s-Lipschitz} and taking into account the choice of
$\sigma$ given by  \emph{ii.}, it follows that:
\begin{equation*}\begin{aligned}
\abs{J_2} &\leq \frac{m}{\sqrt{4\pi}} \int_0^t \frac{|s(t)-s(\tau)|}{ (t-\tau)^{3/2}}\exp\left\{-\frac{|s(t)-s(\tau)|^2}{4(t-\tau)}\right\} d\tau\\
 &\leq \frac{m(|b_0|+2m|b|)}{\sqrt{4\pi}} \int_0^t \frac{1}{ (t-\tau)^{1/2}}d\tau
= \frac{2m(|b_0|+2m|b|)}{\sqrt{4\pi}} t^{1/2}\leq
\frac{2m(|b_0|+2m|b|)}{\sqrt{4\pi}} \sigma^{1/2}\leq \frac{1}{2}.
\end{aligned}\end{equation*}
Before we consider $J_3$, we need the following auxiliary estimates. The inequality $y e^{-y^2}\leq e^{-\frac{y^2}{2}}$
implies that
\begin{equation}\label{formula30}
|G_x(x,t,\xi,\tau)|\leq \frac{1}{\sqrt{4\pi}(t-\tau)}
\exp\left\{-\frac{|x-\xi|^2}{8(t-\tau)}\right\}.
\end{equation}
The definitions of $s(t)$ and
$s_1(\tau)=s(\tau)+v_R\alpha^{-1}(\tau)$, using that $b<0$, and by  the condition \emph{iii},
yield
\begin{equation}\label{case-b-negativo}
|s(t) -s_1(\tau) | \geq \abs{v_R}-|b_0|\sigma>0.
\end{equation}
If we integrate
\eqref{formula30} we get
\begin{equation}\label{estimate-derivative-G}
\begin{aligned}
\int_0^t |G_x(s(t),t,s_1(\tau),\tau)|\,d\tau &
\leq \frac{1}{\sqrt{4\pi}}\int_0^t \frac{1}{t-\tau}\exp\left\{-\frac{|s(t)-s_1(\tau)|^2}{8(t-\tau)}\right\} \,d\tau\\
&\leq \frac{1}{\sqrt{4\pi}} \int_0^t \frac{1}{t-\tau} \exp\left\{-\frac{(\abs{v_R}-|b_0|\sigma)^2}{8(t-\tau)}\right\}\,d\tau\\
&=\frac{1}{\sqrt \pi}
\int_{\frac{\abs{v_R}-|b_0|\sigma}{\sqrt{8t}}}^\infty \frac{1}{z}
e^{-z^2}\,dz \leq\frac{1}{\sqrt \pi}
\int_{\frac{\abs{v_R}-|b_0|\sigma}{\sqrt{8\sigma}}}^\infty
\frac{1}{z} e^{-z^2}\,dz,
\end{aligned}
\end{equation}
where we used the change of variables
$z=\frac{\abs{v_R}-|b_0|\sigma}{\sqrt{8(t-\tau)}}$.  By the last estimate and by condition   \emph{iv}
\begin{equation}\label{J3i}
\abs{J_3}\leq 2m \int_0^t \abs{G_x(s(t),t,s_1(\tau),\tau)} d\tau
\leq \frac{2m}{\sqrt{\pi}}
\int_{\frac{\abs{v_R}-|b_0|\sigma}{\sqrt{8 \sigma}}}^\infty
\frac{1}{z}e^{-z^2}dz\leq \frac{1}{2}.
\end{equation}
 The estimates for $J_i,
i=1,2,3$ establish that $ T(M) \in  C_{\sigma,m}$ since
\[
\|T(M)\| \leq  J_1+J_2+J_3 \leq m, \quad \forall M\in C_{\sigma,m},
\]
by the choice of $m$ in \eqref{choice-m}.\\

It remains to consider the case $b>0$. It is clear that the only
modification needed is the estimte \eqref{case-b-negativo}. For this use
\begin{equation}\label{keyb}
|s(t)-s_1(\tau)|=|s(t)-s(\tau)-v_R\alpha^{-1}(\tau)|\geq
\abs{|v_R|\alpha^{-1}(\tau)-|s(t)-s(\tau)|}\geq
|v_R|-(|b_0|+m)\sigma,
\end{equation}
which may be estimated from below by a positive constant for some
$\sigma$ small enough. Then, we have the same result as in the
case $b<0$ assuming analogous conditions to \emph{i., ii.,
iii., iv.} above. The main difference between the cases $b leq 0$ and $b>0$  is that in  the case $b>0$   in
\eqref{case-b-negativo}  the difference between the free
boundary $s(t)$ and the source $s_1(\tau)$ for $0\leq \tau \leq t$ now
depends  on the bound of
the initial data \eqref{choice-m}. \\

{\it Step 2.-} The mapping $T:C_{\sigma,m}\to C_{\sigma,m}$
defined in \eqref{T} is a contraction for $\sigma$ small enough.
In the sequel the constant $C$ is arbitrary and may change from
line to line. Let $M, \tilde{M} \in C_{\sigma,m}$, and
\begin{align}\label{formula1}
s(t)=s_I-b_0
\left(\sqrt{1+2t}-1\right)-b\int_0^t M(\tau)  \alpha^{-1}(\tau)\,d\tau, \,.\\
\tilde{s}(t)=s_I-b_0
\left(\sqrt{1+2t}-1\right)-b\int_0^t \tilde M(\tau) \alpha^{-1}(\tau)\,d\tau \notag\,.
\end{align}

The following auxiliary estimate holds:
\begin{equation}\label{difference-s}\begin{aligned}
\abs{s(t)-\tilde s(t)}&\leq \abs{b}\int_0^t |M(\tau)-\tilde M(\tau)| \alpha^{-1}(\tau)\,d\tau \leq \abs{b}\|M-\tilde M\| \int_0^t \sqrt{2\tau+1}d\tau\\
&=\frac{\abs{b}}{3} \|M-\tilde M\| \left[(2 t+1)^{3/2}-1\right] .
\end{aligned}\end{equation}
It is straightforward from \eqref{formula1} that
\begin{equation}\label{diference-dot-s}
|\dot s(t)-\dot{\tilde s}(t)|\leq 2\abs{b} \|M-\tilde M\|, \quad 0<t\leq\sigma<1.
\end{equation}
By condition $\emph{i.}$ on $\sigma$ and \eqref{s-Lipschitz} it
follows that
\begin{equation} \label{3i}
\max\{|s(t)-s(\tau) |,|\tilde s(t)-\tilde s(\tau) |\} \leq
(|b_0|+2 m\abs{b})|t-\tau|\leq (|b_0|+2 \abs{b})m|t-\tau|\, .
\end{equation}
To show that $T$ is a contraction we proceed as follows. \begin{align*}
|T(M)-T(\tilde{M})| \leq & \,2\left[ \int_{-\infty}^{s(0)} |u_I'(\xi)| |G(s(t),t,\xi,0)-G(\tilde{s}(t),t,\xi,0)| \,d\xi\right]\\
&+2\left|\int_0^tM(\tau)G_x(s(t),t,s(\tau),\tau)-\tilde{M}(\tau)G_x(\tilde{s}(t),t,\tilde{s}(\tau),\tau)\, d\tau\right|\\
&+ 2\left|\int_0^tM(\tau)G_x(s(t),t,s_1(\tau),\tau)-\tilde{M}(\tau)G_x(\tilde{s}(t),t,\tilde{s}_1(\tau),\tau)\, d\tau\right|\\
=&:\,\Aa_1+\Aa_2+\Aa_3.
\end{align*}
Without loss of generality assume that $\tilde{s}(t) > s(t)$. The
mean value theorem applied to the kernel $G(x,t,\xi,0)$ gives for
some $\bar{s} \in [s(t),\tilde{s}(t)]$
\begin{equation}\label{formula3}
|G(s(t),t,\xi,0)-G(\tilde{s}(t),t,\xi,0)| \leq
|G_x(\bar{s},t,\xi,0)|\cdot|s(t)-\tilde{s}(t)|.
\end{equation}
Recall that
$$
|G_x(\bar{s},t,\xi,0)| = \frac{|\bar{s}-\xi|}{2
t}\frac{1}{\sqrt{4\pi t}}\exp\left\{-\frac{|\bar
s-\xi|^2}{4t}\right\}\leq \frac{1}{\sqrt t} \frac{1}{\sqrt{4\pi
t}} \exp\left\{-\frac{|\bar s-\xi|^2}{8t}\right\},
$$
where we have used the relation $y e^{-y^2}\leq e^{-y^2/2}$.
Hence \eqref{formula3} simply reduces to
\[
|G(s(t),t,\xi,0)-G(\tilde{s}(t),t,\xi,0)| \leq \frac{C}{\sqrt{t}} G(\bar{s}(t),2t,\xi,0))|s(t)-\tilde{s}(t)|.
\]
Integrating in $\xi$, together with \eqref{difference-s}  yields
\[
\Aa_1 \leq C|b|\|u_I'\| \|M-\tilde{M}\| \left\{\frac{(1+2t)^{3/2}-1}{t^{1/2}} \right\}.
\]
Since $\displaystyle\lim_{t \to 0} t^{-1/2}((1+2t)^{3/2}-1)=0$,
for $\sigma$ sufficiently small we have $\Aa_1 \leq
\frac{1}{6}\|M-\tilde M\|$. To estimate $\Aa_2$ we proceed as
follows
\begin{equation*}\begin{aligned}
\abs{\Aa_2} \leq &\,2\left|\int_0^t M(\tau)G_x(s(t),t,s(\tau),\tau)-\tilde{M}(\tau)G_x({s}(t),t,{s}(\tau),\tau) d\tau\right|\\
&+2\left|\int_0^t\tilde{M}(\tau)G_x(s(t),t,s(\tau),\tau)-\tilde{M}(\tau)G_x(\tilde{s}(t),t,\tilde{s}(\tau),\tau) d\tau\right|\\
=&\,:\Aa_{21} + \Aa_{22}.
\end{aligned}\end{equation*}
The Lipschitz bound \eqref{s-Lipschitz} for $s$ yields
$$
|G_x(s(t),t,s(\tau),\tau)|\leq \frac{1}{2\sqrt{4\pi}}
\frac{\abs{s(t)-s(\tau)}}{(t-\tau)^{3/2}}\leq \frac{(|b_0|+2
m\abs{b})}{2\sqrt{4\pi}} \frac{1}{(t-\tau)^{1/2}},
$$
and consequently
\begin{align*}
\abs{\Aa_{21}}\leq &\, 2 \|{M-\tilde M}\|\int_0^t |G_x(s(t),t,s(\tau),\tau)|d\tau \\
\leq & \,C\|{M-\tilde M}\|\int_0^t
\frac{1}{(t-\tau)^{1/2}}d\tau\leq Cm\|{M-\tilde M}\|
\sigma^{1/2}\leq \frac{1}{12}\|M-\tilde M\|,
\end{align*}
for $\sigma$ small enough. To estimate $\Aa_{22}$  proceed as
follows:
\begin{equation*}\begin{aligned}
|G_x (s(t),&t,s(\tau),\tau)-G_x(\tilde s(t),t,\tilde s(\tau),\tau)|\\
&=C\left|\frac{s(t)-s(\tau)}{t-\tau} G(s(t),t,s(\tau),\tau)-\frac{\tilde s(t)-\tilde s(\tau)}{t-\tau} G(\tilde s(t),t,\tilde s(\tau),\tau)\right| \\
&\leq C \left| \frac{s(t)-s(\tau)}{t-\tau}- \frac{\tilde s(t)-\tilde s(\tau)}{t-\tau} \right|G(s(t),t,s(\tau),\tau)\\
&\quad+C \frac{\tilde s(t)-\tilde s(\tau)}{t-\tau} \left| G(s(t),t,s(\tau),\tau)-G(\tilde s(t),t,\tilde s(\tau),\tau)\right| \\
&=:\mathcal B_1+\mathcal B_2.
\end{aligned}\end{equation*}
In order to estimate $\mathcal B_1$ we use the mean value theorem
\begin{equation}\label{formula4}
\frac{[s(t)-\tilde s(t)]-[s(\tau)-\tilde s(\tau)]}{t-\tau}=\dot
s(\bar \tau)-\dot{\tilde s}(\bar \tau)
\end{equation}
for some $0<\bar \tau<t$. By the previous equality and
\eqref{diference-dot-s} we have
$$
\mathcal B_1\leq \frac{C}{(t-\tau)^{1/2}}|\dot s(\bar \tau)-\dot{\tilde s}(\bar \tau)|\leq \frac{C}{(t-\tau)^{1/2}} \|M-\tilde M\|.
$$
On the other hand, to handle the term $\mathcal B_2$, we first
note that
\begin{equation}\label{formula10}\begin{aligned}
&\left| G(s(t),t,s(\tau),\tau)  -G(\tilde s(t),t,\tilde s(\tau),\tau)\right| \\
&\leq G(s(t),t,s(\tau),\tau) \left| 1-\exp\left\{\frac{-(\tilde s(t)-\tilde s(\tau))^2+ (s(t)-s(\tau))^2}{4(t-\tau)} \right\}\right|.
\end{aligned}\end{equation}
Define now
\begin{equation}\label{S}
S:=(s(t)-s(\tau))^2-(\tilde s(t)-\tilde s(\tau))^2 =
\left[s(t)-s(\tau)+\tilde s(t)-\tilde s(\tau)\right][s(t)-\tilde
s(t)- (s(\tau)-\tilde s(\tau))].
\end{equation}
The mean value theorem \eqref{formula4} and the estimate
\eqref{diference-dot-s} lead to
\begin{equation}\label{formula14}\begin{aligned}
\abs{[s(t)-\tilde s(t)]-[s(\tau)-\tilde s(\tau)]}&=\abs{\dot
s(\bar \tau)-\dot\tilde s(\bar \tau)}(t-\tau)\leq C\|M-\tilde
M\|(t-\tau).
\end{aligned}\end{equation}
On the other hand, we recall again the Lipschitz estimate
\eqref{3i}, i.e.,
\begin{equation}\label{equation12}
\max\{|s(t)-s(\tau)|,|\tilde s(t)-\tilde s(\tau)|\}\leq C
m(t-\tau),
\end{equation}
for a constant depending on $|b|$, $|b_0|$, which yields an estimate for \eqref{S},
$$
\frac{\abs{S}}{t-\tau}\leq Cm\sigma \|M-\tilde M\|.
$$
The combination of the above inequality with \eqref{formula10}
together with the mean value theorem shows that
$$
\left| G(s(t),t,s(\tau),\tau)  -G(\tilde s(t),t,\tilde s(\tau),\tau)\right|\leq G(s(t),t,s(\tau),\tau) Cm \sigma\|M-\tilde M\|,
$$
and thus the term $\mathcal B_2$ is estimated using \eqref{equation12}
$$
\mathcal B_2\leq C m^2\|M-\tilde M\| \sigma\frac{1}{(t-\tau)^{1/2}}.
$$
Multiplying $\mathcal B_1+\mathcal B_2$ by $\tilde{M}(\tau)$  and
integrating over the interval $[0,t]$ yields
\begin{equation*}\begin{aligned}
\Aa_{22}&\leq C m\int_0^t |G_x (s(t),t,s(\tau),\tau)-G_x(\tilde s(t),t,\tilde s(\tau),\tau)|\,d\tau\\
&\leq  C m\int_0^t (\mathcal B_1+\mathcal B_2) \,d\tau\leq C m^3\|M-\tilde M\|\sigma^{1/2}<\frac{1}{12}\|M-\tilde M\|,
\end{aligned}\end{equation*}
for $\sigma$ small enough.

The next step is to estimate $\Aa_{3}$. Split the integral into two terms
\begin{equation*}\begin{aligned}
\abs{\Aa_3} \leq &\,2\left|\int_0^t M(\tau)G_x(s(t),t,s_1(\tau),\tau)-\tilde{M}(\tau)G_x({s}(t),t,{s_1}(\tau),\tau) d\tau\right|\\
&+2\left|\int_0^t\tilde{M}(\tau)G_x(s(t),t,s_1(\tau),\tau)-\tilde{M}(\tau)G_x(\tilde{s}(t),t,\tilde{s}_1(\tau),\tau) d\tau\right|\\
=&\,:\Aa_{31} + \Aa_{32}.
\end{aligned}\end{equation*}
The estimate for $\Aa_{31}$ is very similar to that of $J_3$ from
\eqref{J3i}. Indeed, we can analogously deduce
\begin{equation*}\begin{aligned}
\abs{\Aa_{31}}&\leq 2 \|{M-\tilde M}\|\int_0^t |G_x(s(t),t,s_1(\tau),\tau)|d\tau \\
&\leq C \|{M-\tilde M}\|
\int_{\frac{\Lambda}{\sqrt{8\sigma}}}^\infty
\frac{1}{z}e^{-z^2}dz<\frac{1}{12}\|{M-\tilde M}\|,
\end{aligned}\end{equation*}
where we used that $\tilde s(t)-\tilde s_1(\tau)\geq \Lambda>0$
for $\sigma$ sufficiently small with
$$
\Lambda:=
\begin{cases}
|v_R|-|b_0|\sigma & \mbox{for } b<0\\
|v_R|-(|b_0|+m)\sigma & \mbox{for } b>0
\end{cases} .
$$
To bound $\Aa_{32}$ we split
\begin{equation*}\begin{aligned}
|G_x (s(t),&t,s_1(\tau),\tau)-G_x(\tilde s(t),t,\tilde s_1(\tau),\tau)|\\
&=C\left|\frac{s(t)-s_1(\tau)}{t-\tau} G(s(t),t,s_1(\tau),\tau)-\frac{\tilde s(t)-\tilde s_1(\tau)}{t-\tau} G(\tilde s(t),t,\tilde s_1(\tau),\tau)\right| \\
&\leq C \left| \frac{s(t)-s_1(\tau)}{t-\tau}- \frac{\tilde s(t)-\tilde s_1(\tau)}{t-\tau} \right|G(s(t),t,s_1(\tau),\tau)\\
&\quad +C \frac{\tilde s(t)-\tilde s_1(\tau)}{t-\tau} \left| G(s(t),t,s_1(\tau),\tau)-G(\tilde s(t),t,\tilde s_1(\tau),\tau)\right| \\
&=:\mathcal B'_1+\mathcal B'_2.
\end{aligned}\end{equation*}
We observe that $\mathcal B_1'$ is estimated exactly the same way
as $\mathcal B_1$. This is a consequence of
\begin{equation*} 
[s(t)-\tilde s(t)]-[s_1(\tau)-\tilde s_1(\tau)]=[s(t)-\tilde
s(t)]-[s(\tau)-\tilde s(\tau)].
\end{equation*}
We can continue from \eqref{formula4} as before to obtain
\begin{equation}\label{formulaB1'}
\mathcal |B_1'|\leq \frac{C}{(t-\tau)^{1/2}}\|M-\tilde
M\|.
\end{equation}
The estimate for $\mathcal B_2'$ is slightly more involved. We
write
\begin{equation}\label{formula5}
\mathcal B_2'=C\frac{\tilde s(t)-\tilde
s_1(\tau)}{t-\tau}|G(\tilde s(t),t,\tilde s_1(\tau),\tau)|\left(
1-\exp\left\{\frac{S'}{4(t-\tau)}\right\}\right)
\end{equation}
for
\begin{equation*}\begin{aligned}
S':&=-(s(t)-s_1(\tau))^2+(\tilde s(t)-\tilde s_1(\tau))^2 \\
&= \left[\tilde s(t)-\tilde s_1(\tau)+ s(t)- s_1(\tau)\right][\tilde s(t)-\tilde s_1(\tau)- s(t)+ s_1(\tau)]. \\
\end{aligned}\end{equation*}
By the definitions of $s_1$ and $\tilde s_1$ (see \eqref{stefan})
we have that
\begin{equation}\label{formula22}
|\tilde s(t)-\tilde s_1(\tau)-  s(t)+ s_1(\tau)|=|\tilde
s(t)-\tilde s(\tau)- s(t)+ s(\tau)|\leq C\|M-\tilde M\|(t-\tau),
\end{equation}
where in the last inequality we have  used estimate
\eqref{formula14}. On the other hand,
\begin{equation}\label{formula2}
\left|[\tilde s(t)-\tilde s_1(\tau)+ s(t)- s_1(\tau)]\right| \leq
|s(t)-s(\tau)|+|\tilde s(t)-\tilde
s(\tau)|+2|v_R|\sqrt{2\tau+1}\leq Cm\sigma,
\end{equation}
if we use again the Lipschitz estimate \eqref{equation12}. Hence
putting together \eqref{formula22} and \eqref{formula2} we get
again that
$$
\frac{\abs{S'}}{t-\tau}\leq Cm\sigma \|M-\tilde M\|,
$$
and consequently, \eqref{formula5} reduces to
\begin{equation*}
|\mathcal B_2'|\leq Cm \frac{\tilde s(t)-\tilde s_1(\tau)}{t-\tau}
G(\tilde s(t),t,\tilde s_1(\tau),\tau) \sigma \|M-\tilde M\|.
\end{equation*}
Integrating the previous expression, using the inequality $y
e^{-y^2}\leq e^{-y^2/2}$, and noting that $\tilde s(t)-\tilde
s_1(\tau)\geq \Lambda>0$, we can give a very rough estimate that
is enough to our purposes:
\begin{equation}\label{formula20}
\int_0^t |\mathcal B_2'|\,d\tau \leq Cm\sigma \|M-\tilde M\|.
\end{equation}
Thus, from the estimates for $\mathcal B_1'$ and $\mathcal B_2'$
from \eqref{formulaB1'} and \eqref{formula20} respectively,
$$
\abs{\Aa_{32}}\leq Cm \int_0^t (\mathcal B_1'+\mathcal B_2')\,
d\tau\leq C \|M-\tilde M\|\left(m
\sigma^{1/2}+m^2\sigma\right)<\frac{1}{12}\|M-\tilde M\|,
$$
for some suitable $\sigma$ small enough. Then, adding  the
estimates obtained for  $\Aa_1,\Aa_2$ and $ \Aa_3$ yields that $T$
is a contraction satisfying for some $\sigma$ small enough
inversely proportional to $m$:
$$
\|TM-T\tilde M\| \leq \frac{1}{2} \|M-\tilde M\|.
$$
This concludes the proof of Theorem \ref{thm:local-M} as desired.
\end{proof}

\subsection{Recovery of $u$}

Theorem \ref{thm:local-M} shows that we have short time existence
of a mild solution for problem \eqref{stefan} (i.e., a solution in
the integral sense). However, one can easily show that:

\begin{cor}
There exists a unique solution of problem \eqref{stefan} in the
sense of Definition {\rm\ref{defi-solution}} for $t\in[0,T]$.
\end{cor}

\begin{proof}
Once $M$ is known, one can construct  $u$ from Duhamel's formula
\eqref{Duhamel}. The smoothness and decay of $u$ follow
immediately from here. One needs to check also that $u$ has well
defined side derivatives at $s_1$. But this follows from formula
\eqref{little-lemma} applied to $s_1(t)$ and the estimate for
$|G_x(s_1(t),t,s_1(\tau),\tau)|$ that follows similarly as the
calculation in \eqref{estimate-derivative-G}.
\end{proof}

\noindent This completes the proof of Theorem
\ref{short-existence}.


\section{Proofs of the Main Results}

From the previous arguments, see \eqref{choice-m}, it is clear
that the obstacle for long time existence in this case is the
possible blow up in time of $\norm{u_x(\cdot,t)}_\infty$ particularly at
the free boundary, i.e. the blow up of $M(t)$. We now formalize
this idea by showing that we can extend the solution as long as
the firing rate $M(t)$ is bounded.

\begin{prop}\label{prop-derivative-bound}
Suppose that the hypotheses of Theorem {\rm\ref{short-existence}}
hold and that $(u(t),s(t))$ is a solution to \eqref{stefan} in the
time interval $[0,T]$. Assume, in addition, that
$$
U_0:=\sup _{x\in(-\infty,s(t_0-\varepsilon)]}  |u_x(x,t_0-\varepsilon)|<\infty
\qquad \mbox{ and that } \qquad M^*=\sup_{t\in(t_0-\varepsilon,t_0)}
M(t)<\infty\, ,
$$
for some $0<\varepsilon<t_0\leq T$. Then
$$
\sup \left\{ |u_x(x,t)| \mbox{ with } x\in(-\infty,s(t)]\, , \,
t\in[t_0-\varepsilon,t_0) \right\}<\infty\, ,
$$
with a bound depending only on the quantities $M^*$ and $U_0$.
\end{prop}

\begin{proof} Differentiating \eqref{Duhamel} in $x$ yields
\begin{equation*}\begin{aligned}
u_x(x,t) = & \int_{-\infty}^{s(t_0-\varepsilon)} G(x,t,\xi,t_0-\varepsilon) u_x(\xi,t_0-\varepsilon)  d\xi \\
 &-\int_{t_0-\varepsilon}^t M(\tau) G_x(x,t,s(\tau),\tau) d\tau + \int_{t_0-\varepsilon}^t M(\tau)G_x(x,t,s_1(\tau),\tau )d\tau\\
 =: & I_1-I_2+I_3.
 \end{aligned}\end{equation*}
The estimate for $I_1$ is straightforward from heat kernel
properties and it depends only on  $U_0$. Let us deal with $I_2$. Since
$M$ is uniformly bounded in the whole interval $t_0-\varepsilon<t<t_0$, we
get
\begin{equation}\label{I2}
|I_2|\leq C \int_{t_0-\varepsilon}^t |G_x(x,t,s(\tau),\tau)|d\tau.\end{equation}
Next, it is shown in \cite[Eq. (1.16), pag.
219]{Friedman:book} that for any Lipschitz continuous function  $s$,
$$
\int_{t-\varepsilon}^{t} \frac{\abs{x-s(\tau)}}{(t-\tau)}G(x,t,s(\tau),\tau)d\tau\leq C, \quad  t\in(t_0-\varepsilon,t_0),
$$
for some $C$ depending on the Lipschitz constant of $s$, $t_0$ and
$\varepsilon$. Then this formula allows to control the expression in \eqref{I2} in order to bound the term $I_2$ for $t\in(t_0-\varepsilon,t_0)$. However, this bound may depend
on $t_0$ and $M^*$ since the Lipschitz constant of $s$ does, see
\eqref{formula1}.

Finally, the same argument works for $I_3$,  replacing
$s$ by $s_1$ in the previous calculations.
\end{proof}

With this result in hand,  our solutions  can be extended to a
maximal time of existence and, we can characterize this maximal time. The
following result holds no matter the sign of $b$.

\begin{thm}\label{charac}
Suppose that the hypotheses of Theorem {\rm\ref{short-existence}}
hold. Then the solution $u$ can be extended up to a maximal time
$0<T^*\leq \infty$ given by
$$
T^*=\sup\{t>0\,:\, M(t)<\infty\}\,.
$$
\end{thm}

\begin{proof}
Assume that the maximal time of existence of a classical solution
$(u(t),s(t))$ to \eqref{stefan} in the sense of Definition
\ref{defi-solution} is $T^*<\infty$, if not there is nothing to
show. By definition we have $T^*\leq\sup\{t>0\,:\, M(t)<\infty\}$.
Let us show the equality by contradiction. Let us assume that
$T^*<\sup\{t>0\,:\, M(t)<\infty\}$ and then, there exists
$0<\varepsilon<T^*$ such that
$$
M^*=\sup_{t\in(T^*-\varepsilon,T^*)} M(t)<\infty\, .
$$
Let $U_0$ be defined as in Proposition \ref{prop-derivative-bound}
with $t_0=T^*$. Applying Proposition \ref{prop-derivative-bound},
we deduce that $u_x(x,t)$ is also uniformly bounded for $x\in
(-\infty,s(t)]$ and $t\in [T^*-\varepsilon,T^*)$ by a constant,
denoted $U^*$. The same Proposition tells us that $U^*$ only
depends on $M^*$ and also on $U_0$, i.e., the uniform bound of
$u_x(x,T^*-\varepsilon)$ for $x\in (-\infty,s(T^*-\varepsilon)]$.
Therefore, we can now restart by using the local in time existence
Theorem \ref{short-existence} using as initial time any $t_0\in
[T^*-\varepsilon,T^*)$ for a time interval whose length does only
depend on $U^*$. Thus, we can extend the solution $(u(t),s(t))$ to
\eqref{stefan} after $T^*$ and find a continuous extension of
$M(t)$ past $T^*$. We have reached a contradiction, hence the
conclusion of the Theorem follows.\end{proof}

We now show, following Friedman's ideas
\cite{Friedman:book}, that it is possible to extend the solution
for a short (but uniform) time $\varepsilon$ for $b<0$.

\begin{prop} \label{prop-global-bound}
Suppose that the hypotheses of Theorem {\rm\ref{short-existence}}
hold and that $(u(t),s(t))$ is a solution to \eqref{stefan} in the
time interval $[0,t_0)$ for $b<0$. There exists $\varepsilon>0$ small
enough such that, if
\begin{equation}\label{bound-derivative}
\sup_{x\in\mathbb (-\infty,s(t_0-\varepsilon)]} |u_x(x,t_0-\varepsilon)|<\infty,
\end{equation}
for $0<\varepsilon<t_0$ then
$$
\sup_{t_0-\varepsilon <t<t_0} M(t)<\infty.
$$
Although the estimate depends on the bound
\eqref{bound-derivative}, $\varepsilon$ does not depend on $t_0$.
\end{prop}

\begin{proof}
We use the integral formulation  (\ref{M}) for $M$, this time
with initial condition at time $t_0-\varepsilon$ for some fixed $\varepsilon$
chosen below, and $t\in(t_0-\varepsilon,t_0)$. Thus
\begin{equation}\label{M1}\begin{aligned}
 M(t) =& -2\int_{-\infty}^{s(t_0-\varepsilon)}G(s(t),t,\xi,t_0-\varepsilon) u_x(\xi,t_0-\varepsilon)  d\xi \\
 &+2\int_{t_0-\varepsilon}^t M(\tau)G_x(s(t),t,s(\tau),\tau) d\tau - 2\int_{t_0-\varepsilon}^t M(\tau) G_x(s(t),t,s_1(\tau),\tau)  d\tau\\
  =& \; :K_1+K_2+K_3.
\end{aligned}\end{equation}
Since $s(t)\geq s(\tau)$, it follows that
$G_x(s(t),t,s(\tau),\tau)\leq 0$. Moreover, $M\geq 0$, hence
$K_2\leq 0$ and this term can be discarded. To estimate $K_3$ let
$$
\Phi(t):=\sup_{t_0-\varepsilon<\tau<t} M(\tau).
$$
Note that
\begin{equation}\label{K3}
\abs{K_3}\leq \Phi(t)\int_{t_0-\varepsilon}^{t}
\abs{G_x(s(t),t,s_1(\tau),\tau)}d\tau.
\end{equation}
To estimate the derivative  $|G_x(s(t),t,s_1(\tau),\tau)|$ we use
that the nonlinear part of $s$ is an increasing function in the
case $b<0$ as in \eqref{case-b-negativo}. Thus, for $\varepsilon$ small
enough, we conclude that
\begin{equation}\label{key}
s(t)-s_1(\tau)= s(t)-s(\tau)-v_R\alpha^{-1}(\tau)\geq
|v_R|-|b_0|\varepsilon>0\,.
\end{equation}
Hence, we can recall the computations in
\eqref{estimate-derivative-G} to estimate
\begin{equation*}\begin{aligned}
\int_{t_0-\varepsilon}^{t} \abs{G_x(s(t),t,s_1(\tau),\tau)}d\tau \leq C
\int_{\frac{|v_R|-|b_0|\varepsilon}{\sqrt{8(t-t_0+\varepsilon)}}}^\infty\frac{1}{z}e^{-z^2}\,dz\leq
C\int_{\frac{|v_R|-|b_0|\varepsilon}{\sqrt{8\varepsilon}}}^\infty\frac{1}{z}e^{-z^2}\,dz.
\end{aligned}\end{equation*}
It is clear that this last  integral can be made less than $1/2$
for some  small enough $\varepsilon$,   which is independent of  $t_0$.
Substituting the above inequality into \eqref{K3} we have the
estimate $|K_3|\leq \frac{1}{2}\Phi(t)$. Finally, note that
$|K_1|\leq C$ depending on $\sup |u_x(x,t_0-\varepsilon)|$ Combining the estimates for $K_1, K_2 , K_3$ with
\eqref{M1} yields
$$M(t)\leq \frac{1}{2}\Phi(t)+C.$$
Taking the supremum  on the left hand side, we get that $\Phi(t)\leq 2C$,
for all $t\in(t_0-\varepsilon,t_0)$, as desired.
\end{proof}

\begin{remark}
Let us point out that the key estimate \eqref{key} to get the
uniformity of the time interval with respect to $t_0$ comes from
the fact that the nonlinear part of the free boundary $s(t)$ is
monotone increasing in the case $b<0$. For the case $b>0$, instead
of \eqref{key} we got \eqref{keyb}, which makes impossible to get
a uniform estimate with respect to $t_0$ since $m$ will depend on
it.
\end{remark}

\noindent Combining Proposition \ref{prop-global-bound}, Theorem
\ref{charac}, and \ref{prop-derivative-bound} gives global
existence for $b<0$, as summarized in the following result:

\begin{thm}
Let $u_I(x) $ be a non-negative $\mathcal C^1((-\infty,s_I])$
function such that $u_I(s_I)=0$, and $u_I, (u_I)_x$ decay at
$-\infty$. Then there exists a unique  global classical solution
$(u(x,t),s(t))$ of the equation \eqref{stefan} with $b<0$ in the
sense of Definition {\rm\ref{defi-solution}}  with
initial data $u_{I}$. Furthermore, the function $s(t)$ is a
monotone increasing function of $t$ if both $b$ and $b_0$ are
negative.
\end{thm}

With this the proof of our main Theorem \ref{main} is complete. We
emphasize that our Theorem \ref{charac} characterizes the possible
blow-up of classical solutions in finite time as the time of
divergence of the firing rate $N(t)$.


\section{Study of the spectrum}

In this section we study the spectrum of the linear version
$\mu=0$ of \eqref{eq:C}:
\begin{equation*}
p_t-\partial_v (vp)-\partial_{vv} p = N(t) \delta_{v=v_R} \quad\mbox{ on }\, (-\infty,0),
\end{equation*}
where
\begin{equation*}
N(t)=-p_v(0,t), \quad p(0,t)=0.
\end{equation*}

The objective is to solve the eigenvalue  problem
\begin{equation}\label{eig_R}
\left\{\begin{aligned}
&\partial_{vv} p+\partial_v (vp)-p_v(0) \delta_{v=v_R} = \lambda \,p , \quad v\in (-\infty, 0),\\
&p(0)=0,
\end{aligned}\right.
\end{equation}
with eigenfunctions $p(v)$ in the space $L^2_{exp}(\mathbb R)$ defined as
$$
L^2_{exp}(\mathbb R): =\left\{p\in L^2(\mathbb R) \,:\, \norm{p}_{L^2_{exp}(\mathbb R)}<\infty\right\},
$$
with norm
$$
\norm{p}^2_{L^2_{exp}(\mathbb R)} :=\int_{\mathbb R} \lp e^{v^2/2} |p(v)|\rp^2\,dv.
$$
Note that although problem \eqref{eig_R} is only defined in
$(-\infty,0)$, it can be easily extended to $\mathbb R$ by odd
reflection. Following an idea developed in
\cite{price-formation,price-formation2011}, we consider the
equivalent problem to (\ref{eig_R}) defined as
\begin{equation}\label{equation10}
\partial_{vv} p_\lambda+\partial_v (vp_\lambda) = \lambda \,p_\lambda \quad\mbox{in }(-\infty,v_R)\cup (v_R,0),
\end{equation}
with $p_\lambda$ satisfying the following properties:
\begin{enumerate}
\item [\bf{(F1)}]  $p_\lambda \in  L^2_{exp}(\mathbb R)$,

\item [{\bf{(F2)}}] $p_\lambda(0)=0$,

\item [{\bf{(F3)}}] Matching condition: ${p_\lambda}(v_R^+ ) =
{p_\lambda}(v_R^-)$,

\item [{\bf{(F4)}}] Jump condition: $\partial_v{p_\lambda}(v_R^+ )
= \partial_v {p_\lambda}(v_R^- ) +\partial_v p_\lambda(0)$.
\end{enumerate}
We are going to define the solution for \eqref{eig_R} of the form:
\begin{align}\label{gene_sol}
p_\lambda(v)=\chi_{(-\infty,v_R)} p^1(v)+\chi_{(v_R,0)} p^2(v),
\end{align}
where each $p^i(v)$, $i=1,2$, is a linear combination of the two
linearly independent solutions of \eqref{equation10} in $\mathbb
R$, and such that the combination (\ref{gene_sol}) satisfies {\bf
(F1)-(F4)}.

The functions $p^1(v)$ and $p^2(v)$ will be calculated by
a standard  classical method used to compute the spectrum for the
classical Fokker-Planck equation given by $\LL ( p ) = \lambda p$,
$v\in \mathbb R$, with
\begin{align}\label{classicalFP}
\LL( p)  \; := \; \partial_{vv} p+\partial_v (vp).
\end{align}
Define first the eigenspace
$$
L^2_m(\mathbb R)=\left\{p\in L^2(\mathbb R) \,:\, \norm{p}_{L^2_m(\mathbb R)}<\infty\right\},
$$
with norm
$$
\norm{p}^2_{L^2_m(\mathbb R)} :=\int_{\mathbb R} \lp 1+v^2\rp^m |p(v)|^2\,dv.
$$

For completeness we recall a well known result on the spectrum for
the classical operator $\LL$, see for instance
\cite{Gallay-Raugel,Risken}:

\begin{lemma} \label{thm-eigenvalues-FK}
For any $m\geq 0$, the spectrum of the operator $\LL$ defined in
\eqref{classicalFP} on $L^2_m(\mathbb R)$ is given by
$$
\sigma(\LL)=\left\{ \lambda\in\mathbb C \,:\,
\mathfrak{R}(\lambda)\leq \tfrac{1}{2}-m\right\} \cup
\left\{-n\,:\, n\in\mathbb N \cup \{0\}\right\}.
$$
Moreover, if $m>\frac{1}{2}$ and if $n\in\mathbb N\cup \{0\}$
satisfies $n+\frac{1}{2}<m$, then $\lambda_n=-n$ is an isolated
eigenvalue of $\LL$, with multiplicity one, and eigenfunction
given by the $n$-th Hermite polynomial
$$
H_n(v)=(-1)^n e^{v^2/2}\frac{d^n}{dv^n} e^{-v^2/2}.
$$
In particular, the spectrum of the Fokker-Planck operator  $\LL$
in  the space $L^2_{exp}(\mathbb R)$ reduces to the eigenvalues
$$
\lambda=-n,\quad n\in\mathbb N\cup\{0\}.
$$
\end{lemma}
It will be very illustrative to give a sketch of the proof of the
above Lemma in view of the computations for (\ref{eig_R}). \\

\noindent{\emph{Proof of Lemma \ref{thm-eigenvalues-FK}: }}
Given $m\in\mathbb N$, we seek a solution $p\in L^2_m(\R)$ for
\begin{equation}\label{eigenvalues-FK}
p_{vv}+\partial_v(vp)=\lambda p,\quad \lambda \in\mathbb C.
\end{equation}

The Fourier transform of a function $p$ is defined as
\begin{align*}
\hat{p}(\xi) =\frac{1}{\sqrt{2\pi}}  \int_{\R} p(v)e^{-iv\xi}\;dv.
\end{align*}
Then, the Fourier transform of \eqref{eigenvalues-FK} yields the
following first order differential equation
\begin{align*}
- \xi^2 \hat{p} - \xi \,\hat{p}_{\xi} =\lambda \hat{p}, \quad \xi\in \R,
\end{align*}
which has solutions given by
\begin{equation}\label{solution-homogeneous}
\hat{p}(\xi) = \left \{ \begin{array}{ll}
\sqrt{2\pi}\,\alpha \; \xi^{-\lambda}e^{-\xi^2/2} ,\quad \textrm{for} \; \xi > 0 ,\\
\sqrt{2\pi}\,\beta \; (-\xi)^{-\lambda}e^{-\xi^2/2} ,\quad \textrm{for} \; \xi < 0 ,
\end{array} \right.
\end{equation}
for some constants $\alpha,\beta\in \mathbb C$. By inverse Fourier
transform we get:
\begin{equation}\label{solution1}\begin{aligned}
p(v)&=\frac{1}{\sqrt{2\pi}}\int_{\mathbb R} e^{i v \xi} \hat
p(\xi)\,d\xi \\
&=\alpha\int_{0}^{+\infty} e^{i v\xi} \xi^{-\lambda}e^{-\xi^2/2}
\,d\xi +\beta\int_{-\infty}^{0} e^{i
v\xi}(-\xi)^{-\lambda}e^{-\xi^2/2}\,d\xi.
\end{aligned}\end{equation}
This is, for each $\lambda\in\mathbb C$, $p$ is the linear
combination of the two linearly independent solutions of
\eqref{eigenvalues-FK}. The constants $\alpha,\beta,\lambda$ will
be determined from the boundary conditions.

Note that the Fourier transform is an isomorphism from $L^2_m(\R)$
to $H^m(\mathbb R)$. Thus, $\lambda$ belongs to the spectrum of
$\LL$ if and only if the function $\hat p$ from
\eqref{solution-homogeneous} belongs to $H^m(\mathbb R)$. Since
$\hat p$ is sufficiently smooth and rapidly decaying for any $m$,
provided we stay away from the origin, we just need to check if
$\hat p\in H^m(0)$, i.e., all derivatives of $\hat p$ of order
less or equal than $m$ are square integrable near the origin.

Clearly, the values $\lambda=-n$, $n\in\mathbb N\cup\{0\}$ are
special and will be considered later. In order to see if  any
other values of $\lambda \in \C$ are admissible, we compute a
general $n$-derivative of the term $|\xi|^{-\lambda}e^{-\xi^2/2} $
around $\xi=0$ :
$$
\frac{\partial ^{(n)}} {\partial \xi^{(n)}} \left[ |\xi |^{-\lambda}e^{-\xi^2/2}\right] \approx c_{\lambda,n}^{\pm}|\xi|^{-\lambda-n}.
$$
It is easy to see that
\begin{equation*}
\frac{\partial ^{(n)} } { \partial  \xi^{(n)}}  \left[|\xi |^{-\lambda}e^{-\xi^2/2}\right]\in L^2(\R)
\end{equation*}
if and only if $Re(\lambda) < \frac{1}{2} - m$. Since the spectrum of $\LL$ is closed, this shows that
$$
\sigma(\LL) \supset \{\lambda\in\mathbb C \:\, \mathfrak{R}(\lambda)\leq \tfrac{1}{2}-m\},
$$
as claimed.

Now we study the values $\lambda= -n$, $n=1,\ldots,m-1$. Let $p_n$
be the corresponding eigenfunction, with Fourier transform given
by \eqref{solution-homogeneous}. Since $\hat p_n$ belongs to $H^m$
near the origin, its $m$-derivative will be $L^2$ integrable,
while the rest of the derivatives $\hat p_n,\partial_\xi \hat p_n,
\ldots,\partial_{\xi}^{(m-1)}(\hat p_n)$ will be continuous at the
origin. This forces to have a very precise values for
$\alpha,\beta$. In particular, for $n$ even, $\alpha=\beta$, while
for $n$ odd we must have $\alpha=-\beta$. Then we have shown that
$\lambda=-n$, $n=1,\ldots,m-1$ are admissible eigenvalues and the
corresponding eigenspaces are one dimensional, with eigenfunctions
given by the well known Hermite polynomial $H_n$ as
$$
p_n(v)=e^{-v^2/2}H_n(v).
$$
\qed

\

We consider now the original problem \eqref{eig_R} and seek for
solutions $p(v)$ of the form (\ref{gene_sol}). Our first
observation is that the values for $\lambda$ are determined only
by the decay of $p$ as $v\to -\infty$. Consequently, if we impose
that the function $p^1$ belongs to $L^2_{exp}(\mathbb R)$, then
this fixes the possible values of the eigenvalues $\lambda$ as in
Lemma \ref{thm-eigenvalues-FK}.  In particular, there is no continuous spectrum. Moreover, for each
$\lambda_n=-n$, $n\in \mathbb N$, we must have
\begin{align}\label{p1}
 p^1(v)=\alpha H_n(v) e^{-v^2/2},
\end{align}
for some $\alpha\in\mathbb R$.

The difference between our problem and  the classical Fokker-Planck
operator lies in the fact in the interval $(v_R,0)$ all solutions
to the ODE \eqref{eigenvalues-FK} for $\lambda=-n$ as given in
\eqref{solution1} are admissible since the behavior at infinity
does no play any role. Nevertheless, we can find a better
expression for the two linearly independent solutions in
\eqref{solution1}. One of those is the well known (\ref{p1}). The
other solution can be easily found by making the following ansatz:
$$
p^2(v) =e^{-v^2/2} H_n(v) g(v).
$$

By imposing that $p^2(v) $ satisfies (\ref{equation10}) one can
obtain an equation for $g(v)$ that reads as follows:
$$
2H_n'g'-vg'H_n+H_ng''=0.
$$
This equation has the following general solution:
$$
g(v)=\beta_1\int_{v_0}^v \frac{e^{s^2/2}}{H_n^2(s)}\,ds+\beta_2,
$$
for some constants $\beta_1,\beta_2\in\mathbb R$, and where we
have fixed any $v_0\in(v_R,0)$ such that $H_n(v_0)\neq 0$ for the
integral to be well defined. Note that $g$ is well defined for all
$v$ even where the denominator vanishes because the Hermite
polynomials only have single roots. Consequently we define
$$
p^2(v) := \beta_1 e^{-v^2/2}H_n(v) \int_{v_0}^v\frac{e^{s^2/2}}{H_n^2(s)}\,ds +\beta_2e^{-v^2/2}H_n(v) ,
$$
and the eigenfunction corresponding to $\lambda =-n$ is simply
\begin{equation}\label{possible-eigenfunction}
p_n(v)=\left\{\begin{aligned}
\alpha e^{-v^2/2}H_n(v), & \quad v\in(-\infty,v_R), \\
\beta_1 e^{-v^2/2}H_n(v) \int_{v_0}^v\frac{e^{s^2/2}}{H_n^2(s)}\,ds +\beta_2e^{-v^2/2}H_n(v),  & \quad v\in(v_R,0].
\end{aligned}\right.
\end{equation}
for some real constants $\alpha, \beta$
For simplicity, let
$$
\theta_n(v):=H_n(v)\; \int_{v_0}^v \frac{e^{s^2/2}}{H_n^2(s)}\,ds.
$$
It is clear, by doing a careful Taylor expansion, that if $v_1$ is a root of $H_n$, then there exists a finite limit for
$$
\Delta_{v_1,n}:=\lim_{v\to v_1} \theta_n (v)\neq 0.
$$

Now we are ready to check if \eqref{possible-eigenfunction} is an
admissible eigenfunction. In the case $n$ is odd integer, the
Hermite polynomial $H_{2n+1}$ vanishes at zero, but as we have
mentioned,
$$
\theta_{2n+1}(v)\to \Delta_{0,{2n+1}} \neq 0,\quad \textrm{as}\;
v\to 0\quad \textrm{for any}\;n\in\N.
$$
Then in this case condition  {\bf{(F2)}}  is satisfied only when
$\beta_1=0$. Then, if we wish $p_{2n+1}$ to be a continuous
function as stated in condition  {\bf{(F3)}}, we must have
$\alpha=\beta_2$ unless $H_{2n+1}(v_R)=0$ that will be considered
afterwards. The solution constructed this way does not satisfy
condition {\bf{(F4)}}, so we conclude that $2n+1$ is not an
admissible eigenvalue.

On the other hand, let us check if $p_{2n}$ is an admissible
eigenvalue. For even integers it holds that $H_{2n}(0) \neq
0$. Thus we can simply take $v_0=0$. Consequently condition
{\bf{(F2)}} is satisfied if and only if
\begin{equation}\label{f1}
\beta_2=0.
\end{equation}
The matching condition {\bf{(F3)}} implies
\begin{align}\label{f3}
\alpha H_{2n}(v_R) = \beta_1 H_{2n}(v_R) \int_{0}^{v_R}
\frac{e^{s^2/2}}{H_n^2(s)}\,ds.
\end{align}
Here we do need to distinguish two cases: if $v_R$ is not a root
for any $H_{2n}$, then the above equality implies
\begin{align}\label{f3-1}
\alpha = \beta_1 \int_{0}^{v_R} \frac{e^{s^2/2}}{H_n^2(s)}\,ds.
\end{align}
If instead $H_{2n}(v_R)=0$ (note that Hermite polynomials only
have single roots), one can repeat a Taylor expansion around $v_R$
for $\theta_{2n}(v)$ and see that
$$
\theta_{2n}(v) \to {\Delta}_{v_R,2n} \neq 0, \quad \textrm{as}\;
v\to v_R.
$$
Consequently (\ref{f3}) cannot be satisfied for these $n$ such
that $H_{2n}(v_R)=0$.

Using conditions (\ref{f1}) and (\ref{f3-1}) for $p_n$ we get
\begin{equation}\label{eigenfunc2n}
p_{2n}(v)=\beta_1 e^{-v^2/2} H_{2n}(v)\cdot\left\{\begin{aligned}
   \int_{0}^{v_R} \frac{e^{s^2/2}}{H_{2n}^2(s)}\,ds, & \quad v\in(-\infty,v_R), \\
  \int_0^v \frac{e^{s^2/2}}{H_{2n}^2(s)}\,ds,  & \quad v\in(v_R,0].
\end{aligned}\right.
\end{equation}

One can easily check that the jump condition {\bf{(F4)}} is
satisfied if and only if
\begin{align}\label{comp_cond}
H_{2n}(0) = H_{2n}(v_R).
\end{align}

Summarizing, we have shown the following:

\begin{thm}
Consider the operator
\begin{equation*}
\left\{\begin{aligned}
&\partial_{vv} p+\partial_v (vp)-p_v(0) \Delta_{v=v_R} = \lambda \,p, \quad  v\in (-\infty, 0] \\
&p(0)=0,
\end{aligned}\right.
\end{equation*}
subject to conditions {\bf{(F1)}} - {\bf{(F4)}}.
\begin{enumerate}

\item There is no continuous spectrum.

\item The value $\lambda =0$ is an eigenvalue with a
one-dimensional eigenspace spanned by the function
\begin{equation*}
p_\infty(v)=\left\{\begin{aligned}
&e^{-v^2/2} & \quad v\in(-\infty,v_R), \\
&\alpha_0 e^{-v^2/2}\int_v^0 e^{v^2/2}\,dv & \quad v\in(v_R,0],
\end{aligned}\right.
\end{equation*}
for
$$
\alpha_0:=\lp\int_{v_R}^0 e^{v^2/2}\,dv\rp^{-1}.
$$

\item There exists a countable set $S\subset R$ such that for all $v_R\not\in S$, there are no other eigenvalues.

\item If $n$ and $v_R$ happen to satisfy the compatibility
condition \eqref{comp_cond}, then $\lambda = -2n$ is an eigenvalue
with eigenspace of finite dimension spanned by the eigenfunction
$p_{2n}(v)$ defined in \eqref{eigenfunc2n}.
\end{enumerate}
\end{thm}

\begin{remark}
We remark that the steady state $p_\infty (v)$ was previously
obtained in \cite{price-formation,Caceres-Carrillo-Perthame}. In
this last paper, it was also shown exponential decay towards
equilibrium $p_\infty$. However, the speed of convergence is
unknown and the spectral analysis does not seem to give any
insight.
\end{remark}


\noindent\textbf{Acknowledgements.} JAC and MdG are partially
supported by the project MTM2011-27739-C04 DGI-MCI (Spain) and
2009-SGR-345 from AGAUR-Generalitat de Catalunya. MPG is supported
by NSF-DMS 0807636 and DMS-1109682. MS is partially supported by the NSF Grant
DMS-0900909. MS was supported by the sabbatical program of the
MEC-Spain Grant SAB2009-0024.

\end{document}